Title: Floating calculation in Mesopotamia

Author: Christine Proust

# Floating calculation in Mesopotamia

Christine Proust[1]
(CNRS & Université Paris-Diderot, Laboratoire SPHERE & Projet SAW)

**Abstract** – Sophisticated computation methods were developed 4000 years ago in Mesopotamia in the context of scribal schools. The basics of the computation can be detected in clay tablets written by young students educated in these scribal schools. At first glance, the mathematical exercises contained in school tablets seem to be very simple and quite familiar, and therefore, they have little attracted the attention of historians of mathematics. Yet if we look more closely at these modest writings, their simplicity proves deceptive. The feeling of familiarity primarily results from a projection on to the past of thought patterns ingrained from our own early learning. Careful observation of school tablets reveals subtle differences between ancient and modern notions of number, quantity, measurement unit, order, divisibility, algorithm, etc. In fact, we find a completely original mathematical world with rules of its own. This article explores the computation practices developed in scholarly milieus in the light of the mathematical basics that was taught in elementary education.

**Table of content**

[1] The research leading to these results has received funding from the European Research Council under the European Union's Seventh Framework Programme (FP7/2007-2013) / ERC Grant agreement n. 269804. This article is the English translation of Proust, "Du calcul flottant en Mésopotamie", published in *La Gazette des Mathématiciens* 138, p. 23-48 (2014). This English version was presented at the conference "Arith 22", Ecole Normale de Lyon, June 22, 2015. General information on cuneiform mathematics, scribal schools, and place value notation is well known to specialists. However, some of the arguments developed here, especially those concerning the duality between measurement values and numbers in sexagesimal place value notation, are personal and are not the subject of consensus (see more details in [11]).

## 1- Cuneiform mathematics

The history of mathematics in Mesopotamia emerged only recently. In 1930 when François Thureau-Dangin, a French Assyriologist, and Otto Neugebauer, an Austrian-American mathematician, published the first mathematical texts written in cuneiform on clay tablets, the prospects for the history of mathematics changed dramatically. Historians discovered that highly sophisticated mathematics developed in the Ancient Near East over a thousand years before Euclid. Cuneiform mathematics is both the oldest in history and the latest in historiography.

Hundreds of clay tablets covered with cuneiform writing provide an extraordinarily rich source of information about the earliest mathematics which has reached us. The vast majority of these tablets are dated to the Old Babylonian period (ca. 2000-1600 BCE).[2] Small corpora dated to earlier periods (second part of the third millennium) and more recent (late first millennium) are of particular interest to the extent that they shed light on mathematical practices that are completely different from those of the Old Babylonian period.

Generally, the exact origin of the mathematical tablets published in the 1930s and 1940s is unknown. Most of them come from illegal excavations and were bought from antiquities dealers in the early twentieth century by large European and American museums or by private collectors. However, sometimes the mathematical tablets come from legal excavations and are better documented: this is the case, for example, of tablets from Ur, Nippur, Susa, Mari and various sites in the Diyala valley in Northern Mesopotamia. The corpus of known mathematical tablets has increased considerably in recent years with the systematic publication of the numerous school tablets that have accumulated in museum storerooms for decades without arousing any interest among historians. Today, the number of tablets published exceeds 2000.

## 2- The problem of the orders of magnitude

One of the most striking features of this ancient mathematics is the use of sexagesimal place value notation (SPVN hereafter - see Appendix 1c). This numerical system is specific to mathematical texts, but traces of the use of SPVN can be found in other types of texts, for example, in administrative or commercial transactions, where they consist essentially in the remains of graffiti. Sexagesimal place value notation is well known not only to historians of cuneiform mathematics,[3] but also to historians of astronomy. Indeed, the sexagesimal positional notation inherited from Mesopotamia was used in astronomy treatises written in Greek, Latin, Arabic, Syriac, Hebrew, Sanskrit, Chinese, and many other languages, including modern European languages up to recent times. The sexagesimal system has something that is familiar to us since we still use it today in measuring time and angles.

However, this familiarity is deceptive. Behind the apparent similarity between the ancient and modern sexagesimal place value notation, hide subtle but profound differences. To show the importance of these differences, let us simply emphasize two aspects of the ancient sexagesimal place value notation: 1) the position of the units is not specified in the number or, in other word, SPVN is a floating notation; 2) numbers served mainly to perform multiplications and reciprocals. The first of these two aspects is probably the most confusing for the modern reader. Indeed, it is necessary for us to know the place of the units in a number in order to recognize how much this number represents. For example, in our decimal system, when we wish to distinguish "one" from "one thousand" or from "one tenth" we use zeros and a decimal point, and we write respectively "1", "1000" and "0.1". But in the sexagesimal place value notation used in cuneiform texts, the numbers 1, 60 and 1/60 are represented by the same sign, a vertical wedge( ).

[2] In what follows, all dates are BCE.

[3] Sexagesimal place value notation has attracted considerable interest from the beginning of the rediscovery of cuneiform mathematics ([7], [19]).

To be fully aware of this feature of sexagesimal place value notation, it is best to look at a multiplication table found in a scribal school. Figure 1 shows the photographs and a copy of a multiplication table by 9. Looking the left hand column of the table, it is easy to guess that the vertical wedge ( ) represents one, and the angle wedge ( ) represents ten. The right hand column provides the products by 9 of the numbers written in the left hand column. For example, in front of (3), we see (27).

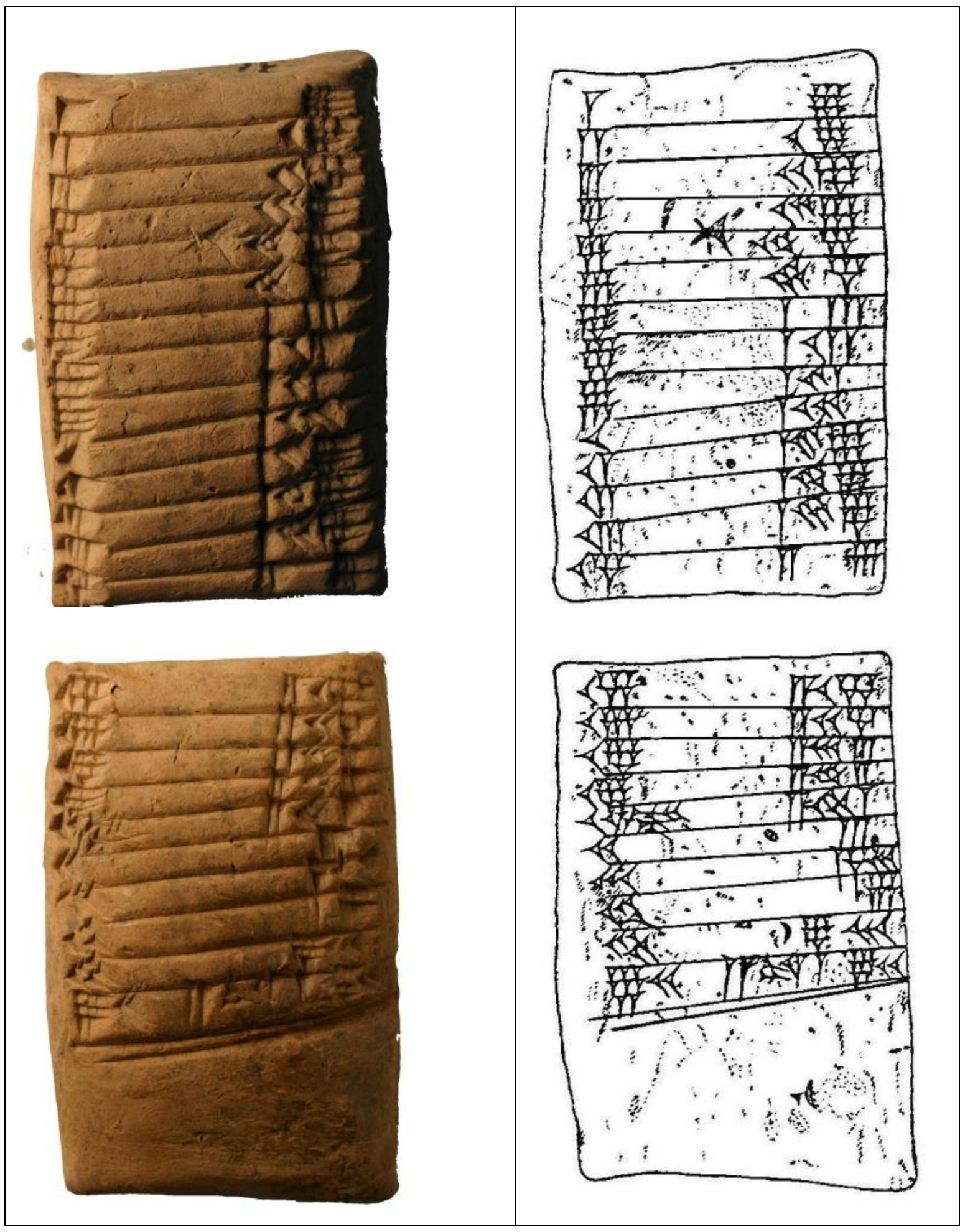

| | |
|---|---|
| 1 | 9 |
| 2 | 18 |
| 3 | 27 |
| 4 | 36 |
| 5 | 45 |
| 6 | 54 |
| 7 | 1:3 |
| 8 | 1:12 |
| 9 | 1:21 |
| 10 | 1:30 |
| 11 | 1:39 |
| 12 | 1:48 |
| 13 | 1:57 |
| 14 | 2:6 |
| 15 | 2:15 |
| 16 | 2:24 |
| 17 | 2:33 |
| 18 | 2:42 |
| 20-1 | 2:51 |
| 20 | 3 |
| 30 | 4:30 |
| 40 | 6 |
| 50 | 7:30 |
| 8:20 | a-ra$_2$ 1 8:20 |

**Figure 1: HS 217a, school tablet from Nippur, multiplication table by 9 (copy Hilprecht 1906, pl. 7, photographs by the author, courtesy Jena University)**

Thus, in front of 7, we expect the number 63. Instead, we see the signs : the sixty is represented by a wedge in the second place, which means that this number is expressed in sexagesimal place value notation. We can transcribe our number 1:3, where the mark ":" is a separator of sexagesimal digits, as in modern digital clocks.[4] The evocation of a digital clock is helpful to understand the ancient sexagesimal system. However, if we continue the reading of the multiplication table, we realize that this comparison has its limits. Indeed, in front of 20, we expect 180, that is, 3×60, or 3:0 in our modern sexagesimal system, but we see 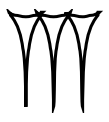. This means that 3×60 is noted with exactly the same sign as the "3" in the third row of the left column. In other words, at the same time the sign 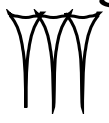 denotes 3 and 3×60. When we see the sign 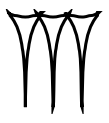 in a multiplication table, we do not know if the wedges represent units, sixties, sixtieth, etc.

Historians were tempted to remedy the uncertainty of the cuneiform notation by adding marks such as comma or zeros (see §5, Table 3, an example extracted from an article by A. Sachs). Which leads to the question: is the lack of graphical systems to determine the place of the units in the number an imperfection in the cuneiform writing system, or on the contrary, is this floating notation an intrinsic property of sexagesimal place value notation?

Consider, for example, the following number, written in SPVN:

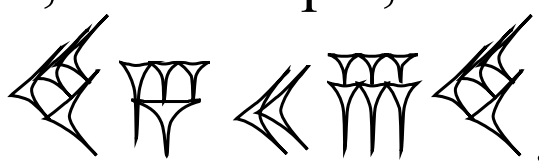.

A modern representation of this number could be: $(44\times60^2+26\times60+40)\times60^p$, $p$ being an indeterminate relative integer. The question could be formulated, in modern language, as follows: should we be concerned about the value of $p$?[5]

To answer this question, we can rely on a particularly rich source of information: the documents produced by the scribal schools. These documents allow us to retrace the main stages of mathematical training of future scribes, and among them, the future "mathematicians". It is therefore possible to have direct access to the mathematical methods that formed the foundation of mathematical erudition. These mathematical basics specifically focused on the concepts related to numbers, measurements values and quantification.[6]

School texts were regarded with some disdain by historians of mathematics because their repetitive exercises appeared to be devoid of attractive mathematical content. School texts have been rehabilitated in recent years as witnesses of the daily life in the scribal schools, that is, as sources for the intellectual history of Mesopotamia, not as sources for the history of mathematics. The collection of tablets preserved in the Istanbul Archaeological Museums shows, on the contrary, the originality of the mathematical concepts taught at Nippur scribal schools.

---

[4] To separate the sexagesimal digits, I chose the separator ':', as in the modern digital documents to make it clear that the base is sixty and not ten. However, the habit in historiography is rather to use the point (44.26.40) or comma (44,26,40) as the separator.

[5] I emphasize that the discussion concerns only SPVN, and not other kinds of notation. For example, numbers used for counting things, say sheep, are not SPVN, but a sexagesiaml non-positional system, where the signs used for 1 and 60 are not the same (see Appendix 1a). For more details and explanations, the reader is invited to read the Appendix containing all relevant information on numerical and metrological systems used in cuneiform mathematics.

[6] I myself have learned cuneiform mathematics by attending a sort of scribal school, since my first steps consisted in discovering, transliterating, translating and commenting hundreds of school tablets from the city of Nippur and kept in the Istanbul Archaeological Museums.

The following sections offer a brief description of the mathematical curriculum in Nippur around the 18th century BCE (see §3), show how the basics were activated to calculate areas (see §4), explain the Babylonian art of floating calculation through the analysis of some numerical algorithms (see §5), and address the resolution of linear and quadratic problems (see §6). The conclusion (see §7) tries to identify the conceptualization of numbers and quantities that emerge from these examples.

## 3- Mathematics taught at Nippur around the 18th century BCE

### The sources

Thousands of clay tablets written by young apprentice scribes inform us about learning mathematics in Mesopotamia. When the archaeological context is known, these tablets were found in places that have often been identified as schools (*edubba* in Sumerian, literally "house of tablets"). Scribal schools may have existed in southern Mesopotamia from the beginning of writing, but the earliest schools left few traces. However, abundant epigraphic sources show that these schools had spread throughout the Middle East at the beginning of the second millennium BC. Many schools were identified in southern Mesopotamia, but the presence of schools is attested over a much larger area, which includes the upper valley of the Euphrates in Iraq, Syria, western Iran and Anatolia.

The mathematical tablets represent, according to the archaeological sites, about 10% to 20% of school tablets. The other school tablets mostly contain evidence of learning cuneiform writing and Sumerian, a language which had not been used as a mother tongue for centuries. While attested over a wide geographical area, the production of these scribal schools is unevenly distributed. The vast majority of known school tablets actually comes from a single site, that of Nippur, where they number in the thousands, including more than 900 mathematical exercises. The Nippur tablets were exhumed by US teams between the late nineteenth century and the first Iraq war, and then divided among several museums (the Istanbul Archaeological Museums, Philadelphia University, the University of Jena, the Baghdad museum and the Oriental Institute of Chicago). School tablets from Nippur form a coherent and quantitatively important set, significant enough to allow statistical consideration. For this reason, the Nippur sources allowed the reconstitution of the curriculum with some precision.[7]In other schools in Mesopotamia, Syria and Iran, significant variations are observed in relation to the shape of the tablets and the layout of the texts they contain. Teaching methods had to be quite diverse. However, it is striking that the content of the most elementary school tablets varies only slightly from site to site. It is likely that the basics taught about the measurement values and numbers in scribal schools were common to all these educational centers.

### The elementary curriculum at Nippur

The school tablets provide information not only through the text they contain, but also through their material appearance and the small notes placed in the margins of the texts (date, scribe's name, catch line, doxology[8] ...). For more details on the methodology that allowed the reconstruction of the curriculum at Nippur, refer to the study by Niek Veldhuis [21]. Applied to mathematics, this method allowed me to reconstruct the general structure of mathematical curriculum (see Table 1).

[7] [2], [10], [12], [16], [17], [20], [21].

[8] The catch line, usually placed at the end of a text, indicates the incipit of the following text in the order of the curriculum. The doxology is a short praise text, usually addressed to Nisaba the goddess of scribes.

| Level | Content |
| --- | --- |
| Elementary | Metrological lists: capacities, weights, surfaces, lengths<br>Metrological tables: capacities, weights, surfaces, lengths, heights<br>Numerical tables: reciprocals, multiplications, squares<br>Tables of square roots and cube roots |
| Intermediate | Exercises: multiplications, reciprocals, surface and volume calculations |

**Table 1: the curriculum at Nippur**

Mathematical education was articulated with literary education (learning cuneiform writing and the Sumerian language). The literary training was to memorize in a fixed order, a series of lists containing, in this order, simple cuneiform signs, Sumerian vocabulary organized thematically then graphically, grammatical structures expressed through paradigmatic Sumerian phrases (the so-called "proverbs"), and finally model contracts used in administration and commerce. The mathematical curriculum included, similarly, memorize lists containing successively capacity, weight, area and length measurement listed in ascending order, tables establishing a correspondence between these measures and sexagesimal positional numbers, and finally numerical tables (inverse, multiplication, square, square roots and cube roots). All these lists and tables providing basic knowledge constitutes what Assyriologists mean by "elementary level" (see Table 1).

After memorizing this set of elementary data, scribes were trained to implement them in some basic numerical algorithms such as multiplication and inversion, and calculate areas and volumes. Subsequently, an "advanced level" would probably start with examples of linear and quadratic problems. However, at least with regard to mathematics, available documentation does not permit us to have a clear idea of the outline of advanced education in Nippur or elsewhere.[9]

As for most of the mathematical texts which look like pure erudition, we can hardly recognize whether they were written by advanced students or by teachers, or to exactly whom they were intended. One of the reasons for this absence of context is that most of the texts that Neugebauer considered "authentically mathematics" (as opposed to the elementary school texts) are of unknown provenance. A few specific insights into advanced education, however, come from texts whose structure seems related to an educational project - among them lists of solved problems commented below (see §6).

A detailed description of the elementary school texts can be found in the different studies mentioned above. Two of these texts deserve special attention in the context of the present story of floating calculation: metrological tables and reciprocal tables.

## Metrological tables

The metrological tables provided young scribes a key mathematical tool to conduct a variety of calculations, including those involving surfaces and volumes. They are presented as tables with two columns. The left column contains a list of capacity, weight, area and length measurements. The whole list of these measurements values, also found in the "metrological lists", provides a complete description of the metrological systems used at the time, not only in mathematics, but also in administrative, commercial or legal activities (see Appendix 2). The right column shows, for each measurement value, a number written in floating sexagesimal place value notation. As an example, Figure 2 shows a copy and translation of a school tablet from Nippur containing an extract of the metrological table for lengths.

[9]The advanced literary education is better documented ([2], [3], [16], [20]).

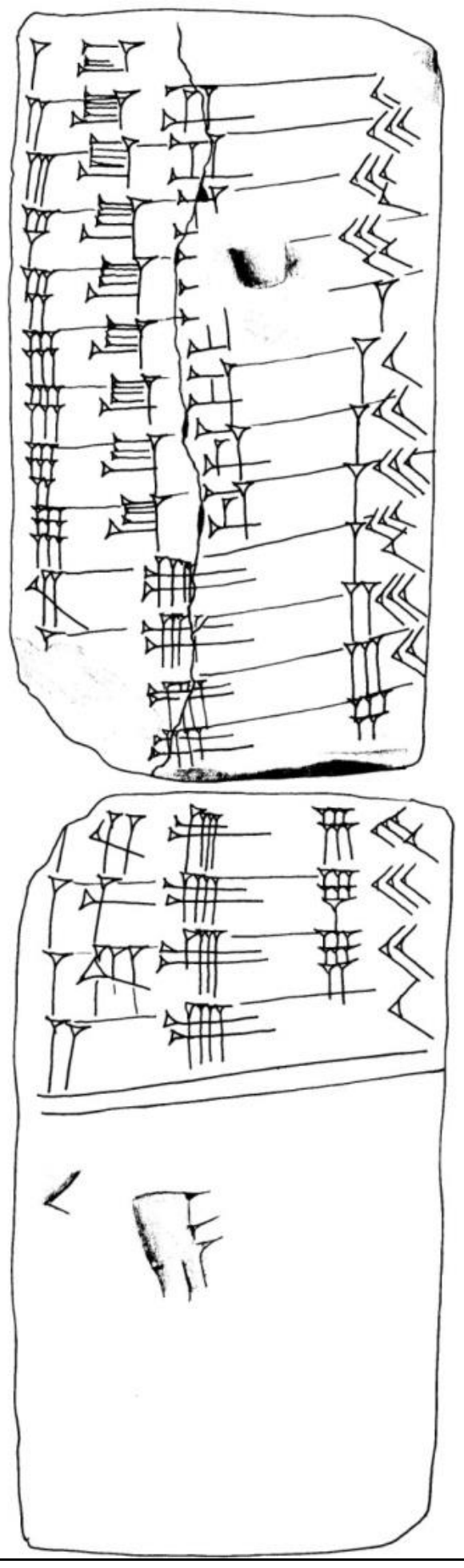

| | |
|---|---|
| 1 *šu-si* | 10 |
| 2 *šu-si* | 20 |
| 3 *šu-si* | 30 |
| 4 *šu-si* | 40 |
| 5 *šu-si* | 50 |
| 6 *šu-si* | 1 |
| 7 *šu-si* | 1:10 |
| 8 *šu-si* | 1:20 |
| 9 *šu-si* | 1:30 |
| 1/3 *kuš* | 1:40 |
| 1/2*kuš* | 2:30 |
| 2/3 *kuš* | 3:20 |
| 5/6 *kuš* | 4:10 |
| 1 *kuš* | 5 |
| | |
| 1 1/3 *kuš* | 6:40 |
| 1 1/2*kuš* | 7:30 |
| 1 2/3 *kuš* | 8:20 |
| 2 *kuš* | 10 |

**Figure 2: HS 241, school tablet from Nippur, extract of the metrological table for length (copy by the author)**

Note the cyclical nature of the right hand column: the number 10 appears in the first line, as a positional number corresponding to 1 *šu-si* (about 1.7 cm.), and in the last line, as a positional number corresponding to 2 *kuš* (about 1 m). Identifying the correct correspondence between a number in SPVN and a measurement value (reading the table from the right column to the left hand column) therefore required a mastery of orders of magnitude. Such a skill was one of the main goals of mathematics education at Nippur. Some examples are described below (see the end of this paragraph and §6).

## Reciprocal tables

The numerical tables studied in the elementary curriculum at Nippur (and other schools) were the following, in this order: reciprocal table, multiplication table by 50, 45, 44:26:40, 40, 36, 30, 25, 24 , 22:30 20, 18, 4:40 p.m., 16, 15, 24:30, 12, 10, 9, 8:20, 8, 7:30, 7:12, 7, 6:40, 6, 5 , 4:30, 4, 3:45, 3:20, 3, 2:30, 2:24, 2, 1:40, 1:30, 1:20, 1:15, squares, square roots, and cube

roots. The first was the reciprocal table. This position ahead of the set of numerical tables would be difficult to explain by pedagogical reasons, and instead seems to reflect the importance of the reciprocals in sexagesimal calculation. Indeed, in mathematical contexts, the division by a number was usually performed by mean of the multiplication by the reciprocal of this number. The set of multiplication tables, which appear to be enumerated in descending order of their main number, could be interpreted instead as a set of "division" tables. For example, the multiplication table by 44:26:40, which is one of the first tables of the set, is also (and perhaps above all) a table of division by 1:21, the reciprocal of 44:26:40.

All the samples of reciprocal tables found in Mesopotamia exhibit more or less the same content: they provide the reciprocals of regular[10] 1-place value numbers, plus those of few 2-place value numbers of frequent use. The reciprocal table consists in a list of clauses: the reciprocal of 2 is 30, the reciprocal of 3 is 20, etc. (see Figure 3 below). The floating nature of the numbers appears here clearly: only the successive sexagesimal digits are given in the table, without any mark indicating the position of the units. In other words, two numbers form a pair of reciprocals if their product is 1, this "1" representing any power of 60

| **Obverse** | | | | **Reverse** | | | |
|---|---|---|---|---|---|---|---|
| The reciprocal of | 2 | is | 30 | The reciprocal of | 48 | is | 1:15 |
| The reciprocal of | 3 | is | 20 | The reciprocal of | 50 | is | 1:12 |
| The reciprocal of | 4 | is | 15 | The reciprocal of | 54 | is | 1:6:40 |
| The reciprocal of | 5 | is | 12 | The reciprocal of | 1 | is | 1 |
| The reciprocal of | 6 | is | 10 | The reciprocal of | 1:4 | is | 56:15 |
| The reciprocal of | 8 | is | 7:30 | The reciprocal of | 1:21 | is | 44:26:40 |
| The reciprocal of | 9 | is | 6:40 | | | | |
| The reciprocal of | 10 | is | 6 | | | | |
| The reciprocal of | 12 | is | 5 | | | | |
| The reciprocal of | 15 | is | 4 | | | | |
| The reciprocal of | 16 | is | 3:45 | | | | |
| The reciprocal of | 18 | is | 3:20 | | | | |
| The reciprocal of | 20 | is | 3 | | | | |
| The reciprocal of | 24 | is | 2:30 | | | | |
| The reciprocal of | 25 | is | 2:24 | | | | |
| The reciprocal of | 27 | is | 2:13:20 | | | | |
| The reciprocal of | 30 | is | 2 | | | | |
| The reciprocal of | 32 | is | 1:52:30 | | | | |
| The reciprocal of | 36 | is | 1:40 | | | | |
| The reciprocal of | 40 | is | 1:30 | | | | |
| The reciprocal of | 45 | is | 1:20 | | | | |

**Figure 3: reciprocal table (MS 3874), copy Friberg ([4], p. 69), and partial translation (first and last item are omitted)**

This list of elementary reciprocal pairs was learnt by heart during the education of the scribes, and the scholars constantly used it in the implementation of certain algorithms (see §5).

[10]A number is regular in base 60 if its reciprocal has a finite number of digits in this base. So the regular numbers in base 60 are those whose decomposition in prime factors contains no factors other than 2, 3 or 5, the prime divisors of 60. The ancient scribes avoided irregular numbers, which were labeled as "without reciprocal" (in Sumerian 'igi nu', where 'igi' means reciprocal, and 'nu' is the negation).

## 4- Computing and quantifying

What are the basic mathematical concepts that emerge from this highly structured set of lists and tables? A first observation touches the operations: in elementary education, there is no mention of addition. Key operations are reciprocals and multiplications; the other operations, square, square root and cube root, derive from multiplication. It is clear here that the basic arithmetical operations are not the classic quartet addition-subtraction-multiplication-division found in modern elementary textbooks, but rather multiplication and reciprocals.

A second observation concerns the numbers: two kinds of numbers are introduced successively by the ancient masters. First, within metrological lists, come the numbers that are used to express measurement values (or numbers of items), that is to say, the numbers expressing quantities. These numbers belong to different non-positional systems. Then, in the metrological tables, a correspondence between measurement values and numbers in sexagesimal place value notation appears. Finally, numerical tables provide a vast repertoire of elementary results of multiplications and reciprocals, in which only numbers in SPVN are involved.[11]

It becomes clear, just by the observation of school material, that the use of floating sexagesimal place value notation is closely linked to multiplication and reciprocals. This link is obviously related to the fact that, in performing a multiplication, the place of the units in a number does not matter, in the sense that it does not affect the list of digits of the result - the only sense that matters in floating notation. So, at this stage, we discover the lineaments of a numerical world of floating numbers on which multiplications and reciprocals act. In this world, the lack of identification of the position of the units in the notation of numbers is not a weakness of the writing system. On the contrary, as we shall see, it is a key feature that enables the implementation of powerful computational algorithms (see §5).

However, with this floating system, how did they count and add sheep? The answer is that sheep are never counted with numbers in SPVN, but only using the additive sexagesimal system, the System S (see Appendix 1a). Sheep are counted and added (or subtracted) with numbers written in an additive system, where the order of magnitude is defined by the form of the signs or other means. These additions of sheep, goats, workers, or various goods, are omnipresent in administrative texts, as well as in some advanced mathematical texts, albeit in a more abstract way. However, additions and subtractions are absent from school exercises: adding sheep is outside the scope of elementary mathematics in scribal schools. These operations on goods appear only in the more professional part of education, for example when the young scribes are trained to write contracts.

To sum up, according to school texts from Nippur, two sorts of numbers fulfill two distinct functions: quantifying and computing. The quantification is ensured by numbers noted in additive systems, and their order of magnitude is well defined. Computation (actually, multiplication and related operations) is ensured by floating numbers written in sexagesimal place value notation. As they are not used for quantification, and their order of magnitude is not defined, these numbers do not denote quantities, but they are just successions of digits, which are never followed by concrete entities such as measurement units or names of items. In this sense, these numbers can be qualified, as Thureau-Dangin did, as "abstract numbers".[12]

[11] The duality between two kinds of numbers, namely non positional numbers for quantification, and positional numbers for multiplication and related operations, is shown in [11].
[12] These ideas are developed and justified with more details and examples in [11].

## Computing surfaces and volumes

Exercises found at Nippur demonstrate how the two kinds of numbers are involved at specific stages of the calculation of surfaces. Figure 4 shows a copy and translation of one of these exercises.

| | |
|---|---|
| 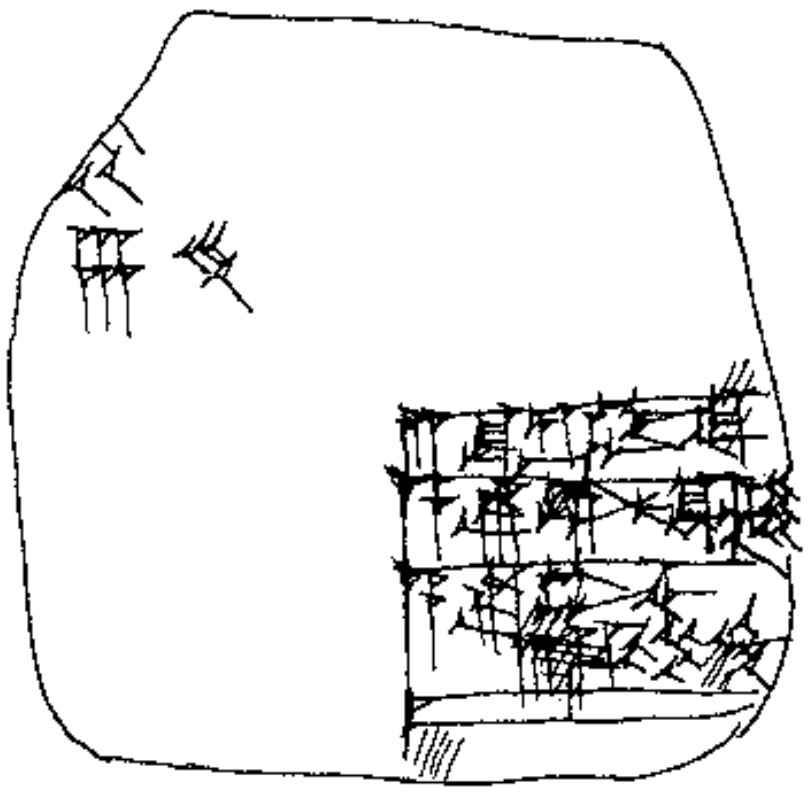 | 20<br>20<br>6:40<br><br>2 *šu-si* is the equal-side.<br>What is its surface?<br>Its surface is 1/3 *še* |

**Figure 4: School tablet from Nippur (UM 29-15-192, copy [9], 248, 251)**

The problem seems to be very simple. In the lower right corner, a short text states that the side of a square ("equal-side") is 2 *šu-si* (about 3 cm) and asks to find its surface, then, it gives the answer, the surface is 1/3 *še*. All quantities are expressed using length or surface measurements, denoted in the same way as in the metrological lists and tables. In the upper left corner, there are three abstract numbers written one under the other. The first two numbers are partially broken, but the visible rests, as well as other similar known tablets, allow us to guess that the third number is the product of the other two; here 20 × 20 = 6:40. According to the metrological table of lengths, the measurement value 2 *šu-si* corresponds to the number 20. This abstract number 20 is multiplied by itself in the upper left corner. The product is 6:40. According to the metrological table for surfaces, used in conjunction with a mental estimate of the order of magnitude, the abstract number 6:40 corresponds to 1/3 *še*.

This exercise, and similar ones, shows that the calculation of surfaces follows three stages:
- Transformation of the measurements of length into abstract number (using the metrological table for lengths)
- Multiplication of two abstract numbers (using multiplication tables)
- Transformation of an abstract number into measurement of surface (using the metrological table for surfaces and an estimation of the order of magnitude).

The calculation of volume and weight followed the same pattern. The system of units of volume used in mathematical texts was shaped artificially, probably in the second half of the third millennium, from the area units: a unit of volume is a unit of area affected by a constant thickness of 1 *kuš* (about. 50 cm). Special metrological table for heights, where 1 was the abstract number associated to 1 *kuš* (instead of 5 in the metrological table for lengths) was used. This trick made it possible to use the same metrological table for surfaces and for volumes.[13]A sophisticated system of coefficients, provided as abstract numbers in coefficient tables, allowed practitioners to calculate the weight of various substances of different densities using metrological tables for weights.

[13] For more explanations, see [10], ch. 6.6.

Without going into more detail, these examples show how the calculation with quantities was based on repeated transformations from measurement values into abstract numbers and vice versa.

## 5- Floating calculation in action: some examples

### Finding the reciprocal of a number

The algorithms that were used for the extraction of reciprocals, square roots and cube roots reveal the potency of floating calculation. I have analyzed these algorithms in detail, focusing on how the scribes controlled the validity of the process, in Karine Chemla's book on *the history of mathematical proof in ancient traditions* ([1], [14]). I will limit myself here to the extraction of reciprocals.

Elementary reciprocals provided by tables (see Figure 3), were known by heart by the experienced scribes. However, how is it possible to calculate the reciprocal of a regular number not given by this table? Again, the answer is to be found in the student exercises, for example the one shown in Figure 5.

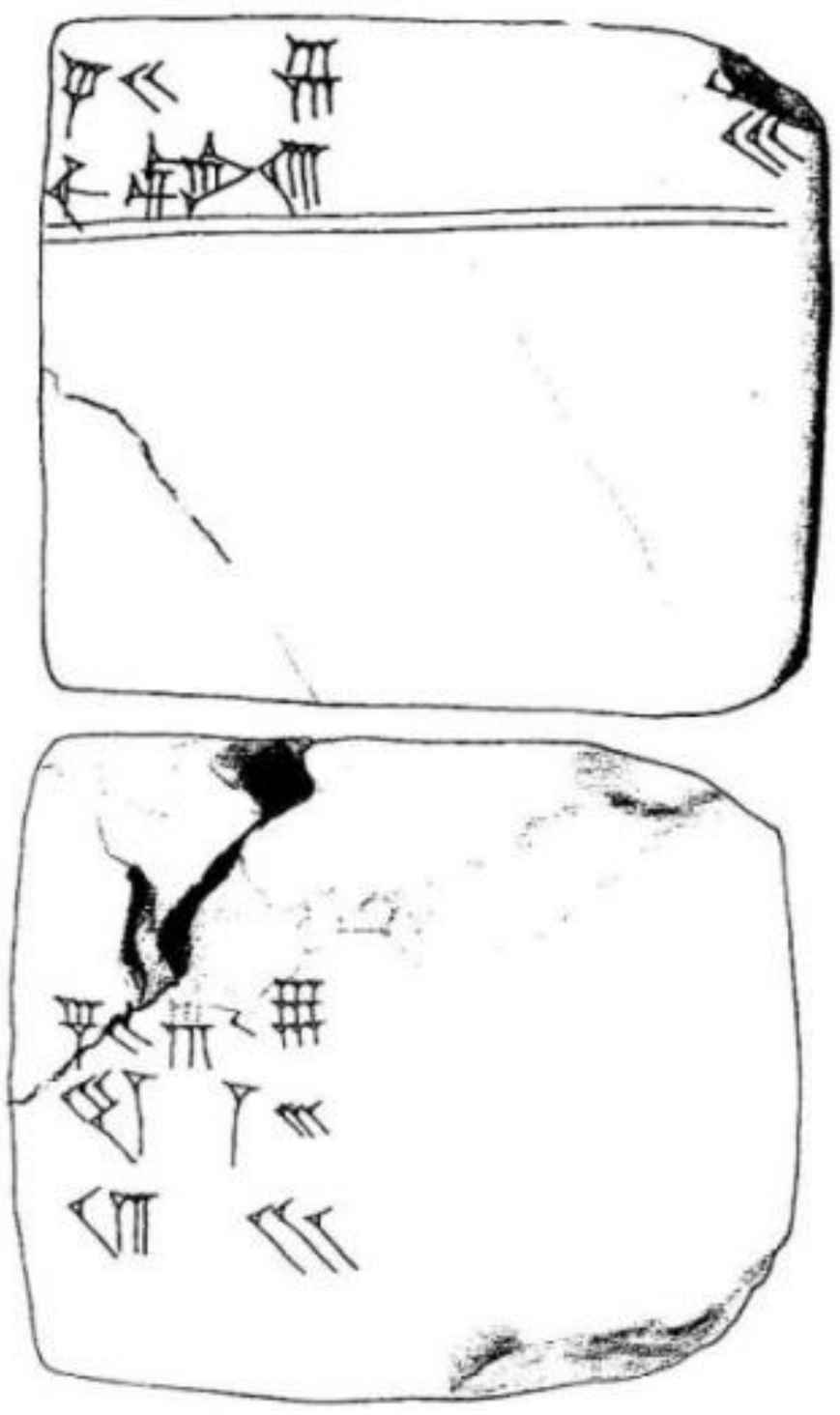

**Obverse**

4:26:40
Its reciprocal 13:30

**Reverse**

4:26:40 9
41[sic] 1:30
13:30

**Figure 5: school tablet from Nippur containing the extraction of a reciprocal (Ist Ni 10241, copy of the author)**

On the obverse, it is stated that the reciprocal of 4:26:40 is 13:30, which is correct (to check, just multiply 1:30 by 4:26:40, the product is 1). On the reverse, the detailed calculation appears, which I reproduce below, correcting the error in the second line (the expected number is 40, while 41 is noted on the tablet):

| | | |
|---|---|---|
| 4:26:40 | | 9 |
| 40 | | 1:30 |
| | 13:30 | |

Here is the explanation of the calculation: 4:26:40, the number for which the reciprocal is sought, ends with the regular number 6:40, so 4:26:40 is "divisible" by 6:40 (I come back to this notion of divisibility below). To divide 4:26:40 by 6:40, we must multiply 4:26:40 by the reciprocal of 6:40. The reciprocal of 6:40 is 9 according to the table of reciprocals (see Figure 3). The number 9 is placed in the right hand sub-column. 4:26:40 multiplied by 9 is 40, thus, 40 is the quotient of 4:26:40 by 6:40. This quotient is placed in the left hand sub-column. The reciprocal of 40 is 1:30 according to the table of reciprocals (see again Figure 3). The number 1:30 is placed in the right hand sub-column. To find the reciprocal of 4:26:40, one only has to multiply the reciprocals of the factors of 4:26:40, that is to say, the numbers 9 and 1:30 placed in the right hand column. This product is 13:30, the reciprocal sought.

To sum up, the left hand sub-column tells us that 4:26:40 can be decomposed into a product of regular factors as follows:

4:26:40 = 6:40 × 40

The right hand sub-column is composed of the reciprocals of these factors, so that the number sought is the product of these reciprocals:

9 × 1:30 = 13:30

The algorithm is based on the decomposition of the number whose reciprocal is sought into regular factors. Thus, the implementation of the algorithm involves the use of some criteria of divisibility.

The notion of divisibility needs to be slightly reworded in the context of the floating calculation. First of all, only divisibility by regular numbers is considered. Another difficulty is that there is no difference between integers and non-integers in the floating system. Any number is divisible by any regular number (for example, 2 can be divided by 5, this gives 24). I take here the word divisible in the following sense: the number *a* is divisible by the regular number *b* if the product of *a* by the reciprocal of *b* is "smaller" than *a*. But what is the meaning of "smaller" in a context where numbers have no order of magnitude? In fact, the goal is to obtain a number which, after division, is simpler than *a*, i.e., a number which has fewer digits than *a* (40 is smaller than 4:26:40), or, if the number of digits is the same, a number less than a (40 is smaller than 50). How can a regular divisor be recognized? Regular divisors appear to form the last portion of the number (for example, 6:40 is a divisor of 4:26:40).[14] We use the same criterion in base 10: the divisibility of a number by 2, 4, 5, 10, 25 etc. appears when looking at the last digits of this number. Sometimes the regular "trailing part" is slightly hidden (see Table 2, where the 3.45 is identified as a factor of 8.45).

When the quotient found is not in the table of reciprocals, the process must be iterated. Such iterations are attested in a tablet kept at Philadelphia University in the United States under the accession number CBS 1215. This tablet is composed of 21 sections containing extractions of reciprocals. The 20th section is reproduced and explained in Table 2 below (I have highlighted in bold regular factors that appear in the final part of the number whose reciprocal is sought, and I underlined the reciprocal found)[15].

[14] Hence the name "trailing part algorithm" given by Friberg to this factorization method.
[15] For a detailed analysis of the text, see [15].

| Transcription | | Explanation | |
|---|---|---|---|
| 5:3:24:2**6:40** | [9] | $n \to$ 5:3:24:2**6:40** | |
| 45:30:**40** | 1:30 | Factors of $n$ | factors of $inv(n)$ |
| 1:8:**16** | 3:45 | 6:40 | 9 |
| 4:**16** | 3:45 | 40 | 1:30 |
| **16** | 3:45 | 16 | 3:45 |
| 14:3:45 | | 16 | 3:45 |
| 5[2:44]:3:45 | | 16 | 3:45 |
| 1:19:6:5:37:30 | | Products | |
| <u>11:51:54:50:37:**30**</u> | 2 | 14:3:45 | |
| 23:43:49:41:**15** | 4 | 5[2:44]:3:45 | |
| 1:34:55:18:**45** | 16 | 1:19:6:5:37:30 | |
| 25:18:**45** | 16 | 11:51:54:50:37:**30** | |
| 6:**45** | 1:20 | | |
| **9** | 6:[40] | $n \to$ 11:51:54:50:37:**30** | |
| 8:53:20 | | **30** | **2** |
| 2:2:22:2:13: 20 | | **15** | 4 |
| 37:55:33:20 | | **3:45** | 16 |
| 2:31:42:13:20 | | **3:45** | 16 |
| 5:3:24:26:40 | | **45** | 1:20 |
| | | **9** | 6:[40] |
| | | Products | |
| | | 8:53:20 | |
| | | 2:22:13: 20 | |
| | | 37:55:33:20 | |
| | | 2:31:42:13:20 | |
| | | 5:3:24:26:40 | |

**Table 2: CBS 1215 #20**

It can be observed that this calculation, besides the iteration of the factorization process, shows a new sophistication. Indeed, once found, the reciprocal is in turn reversed by the same algorithm, which obviously produces the initial number. This loop illustrates the following arithmetical property: the reciprocal of the reciprocal of a number is this number itself.[16]

The simplicity of the algorithm, and therefore its potency, is essentially due to the floating system of notation. This feature of cuneiform sexagesimal place value notation has not always been perceived by historians, as evidenced by the excerpt from the study of CBS 1215 made by Abraham Sachs reproduced in Table 3.

[16] See more details in [15].

| $c$ | $B$ | $a$ | $\bar{a}$ | $1+\bar{a}b$ | $\overline{1+\bar{a}b}$ | $\bar{c}$ |
|---|---|---|---|---|---|---|
| […]<br>#20 | | | | | | |
| 5,3,24,26;40<br>45,30,40<br>1,8,16<br>4,16 | 5,3,24,20<br>45,30, 0<br>1,8,0<br>4,0 | 6;40<br>40<br>16<br>16 | 0;9<br>0;1,30<br>0;3,45<br>0;3,45 | 45,30,40<br>1,8,16<br>4,16<br>16 | 0;0,0,1,19,6,5,37,30<br>0;0,0,52,44,3,45<br>0;0,14,3,45<br>0;3,45 | 0;0,0,0,11,51,54,50,37,30<br>0;0,0,1, 19,6,5,37,30<br>0;0,0,52,44,3,45<br>0;0,14,3,45 |
| 0;0,0,1,19,6,5,37,30<br>23,43,49,41,15<br>1,34,55,18,45<br>25,18,45<br>6,45 | 0;0,0,0,11,51,54,50,37<br>23,43,49,41,0<br>1,34,55,15,0<br>25,15,0<br>6,0 | 0;0,0,0,0,0,0,0,0,30<br>15<br>3,45<br>3,45<br>45 | 2,0,0,0,0,0,0,0,0<br>0;4<br>0;0,16<br>0;0,16<br>0;1,20 | 23,43,49,41,15<br>1,34,55,18,45<br>25,18,45<br>6,45<br>9 | 0;0,0,0,0,2,31,42,13,20<br>0;0,0,0,0,37,55,33,20<br>0;0,0,2,22,13,20<br>0;0,8,53,20<br>0;6,40 | 5,3,24,26;40<br>0;0,0,0,0,2,31,42,13,20<br>0;0,0,0,0,37,55,33,20<br>0;0,0,2,22,13,20<br>0;0,8,53,20 |

**Table 3: The « technique » of extraction of reciprocals in CBS 1215 #20 as explained byA. Sachs ([18], pp. 238-240)**

Sachs'calculations are unnecessarily complicated. First, he tried to express the algorithm by an algebraic formula. Then he considered that the number $c$, whose reciprocal is sought, was broken down into a sum $a + b$, and that the reciprocal was calculated according to the "formula" $\bar{a}(\overline{1+\bar{a}b})$, noting $\bar{a}$ the reciprocal of $a$. Second, he considered the numbers as quantities, to which he restored the order of magnitude by adding a lot of "zeros". This example shows how the use of modern concepts may be inappropriate and hinder a full understanding of old procedures.

## Floating numbers and fixed values: the example of the surface of a disk

The relationship between abstract numbers and measurement values is particularly interesting to observe in elementary geometric situations. As we have seen, calculating a surface is a basic exercise in the curriculum. However, the relationship between abstract numbers and measurement values is not always explicitly shown, as in the text commented before (see §4). Many exercises dealing with sides and surfaces do not mention any measurement values. These exercises often contain diagrams with abstract numbers representing the linear measurements (usually noted outside of the figure) and the corresponding surface (generally noted within the figure). The tablet shown in Figure 6, of unknown origin, exhibits a circle and abstract numbers corresponding to the perimeter (above the circle), the square of the perimeter (right of the circle), and the surface (inside the circle).

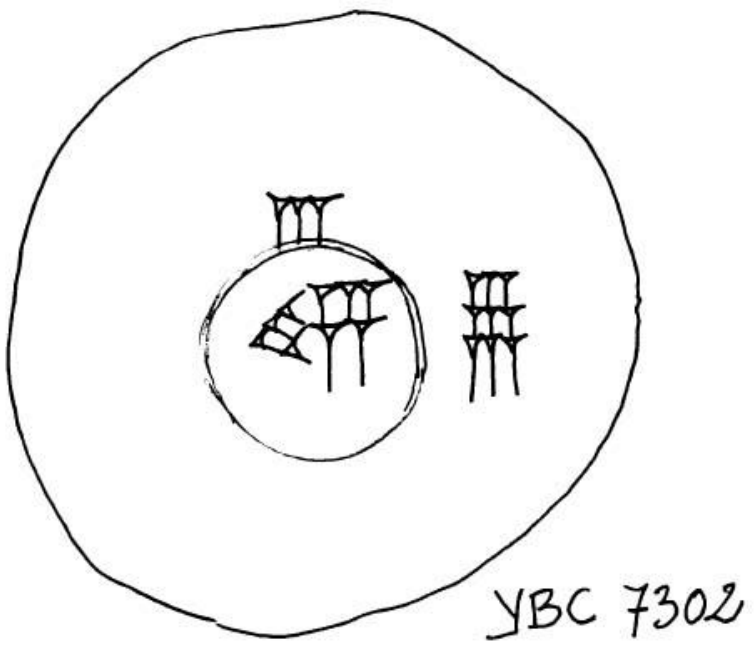

**Figure 6: diagram with numbers (YBC 7302, school tablets, Yale University; copy by the author)**

It appears clearly that the surface (45) is obtained by multiplying the square of the perimeter (9) by 5. This is a general rule, found in many mathematical texts.

Let us explain this relation in modern language. If the area is noted *a,* the perimeter is noted $p$, the relation between $p$ and $a$ is:

$$a = p \ / 4\pi.$$

The scribes roughly estimated the ratio of the perimeter of a circle to the diameter as 3 (in modern language: $\pi \approx 3$). Thus:

$a \approx p / 12$

As the reciprocal of 12 is 5, the area was evaluated by the ancient scribes, in floating notation, as $p \times 5$.

It is notable that, in this tablet, no measurement unit is specified. So, what do the numbers represent? Should we interpret the number 3 as representing a measurement value, such as a perimeter 3 *ninda*? We have seen that, in the metrological tables, the numbers in the right hand sub-column describe several cycles. For example, the number 3 may correspond to several different length measures:

½ *kuš* 3 *šu-si* → 3
3 *ninda* → 3
3 *uš* → 3
6 *danna* → 3

Does the text refer to one of these correspondences, for example 3 *ninda*, excluding the others? This is the interpretation generally found in modern publications. However, no indication allows us to interpret the text in this way. The diagram and numbers may simultaneously represent all the possible metric situations listed in Table 4.

| Perimeter (3) | Surface (45) | Order of magnitude |
|---|---|---|
| | | |
| ½ *kuš* 3 *šu-si* | 2 ¼ *še* | A tablet |
| 3 *ninda* | 2/3 *sar* 5 *gin* | An orchard |
| 3 *uš* | (1×18 + 1×6 +3) *gan* | A city |
| 6 *danna* | (1×1080 + 3×180)*gan* | A region |

**Table 4: possible orders of magnitude corresponding to the data of YBC 7302**

As in the case of the tables of reciprocals, floating notation confers a great flexibility to the text.

In the usual modern interpretations, the number 3 is understood as “3 *ninda*”, which changes the meaning of the text. Indeed, numbers corresponding to lengths are usually interpreted as measurement values where the unit "*ninda*" is omitted, but implied. This view is expressed by Høyrup as follows: “The basic measure of horizontal distance is the nindan ("rod"), equal to c. 6 m. Mostly, this unit is not written but remains implicit.” ([6], p. 17).[17] The adjunction of an “implicit” measuring unit, for example replacing “3” with “3 *ninda*”, may appear as anodyne. However, by adding a measuring unit, the modern historian changes the very nature of the numbers. While 3 is written in SPVN, floating in nature, the number “3” in “3 *ninda*” has a fixed order of magnitude (3 units) and denotes a quantity.

[17] I argued elsewhere [11] against the idea that the abstract numbers can be interpreted as measurement values whose unit was omitted.

## 6- Linear and quadratic problems in YBC 4663

The constant switching from measurement values to abstract numbers and vice versa is quite apparent in elementary school texts, where only multiplications and reciprocals are involved. How are measurement values and abstract numbers articulated in more sophisticated mathematical practices? Once more, pedagogical material provides some elements in answer to this question. Some interesting evidence is offered by tablet YBC 4663, which was probably used at the beginning of advanced education in Southern Mesopotamia during the Old Babylonian period.[18] The tablet contains a list of eight problems with detailed procedures of resolution. To be solved, problems 1 to 6 require only sequences of multiplications and reciprocals. The last two are quadratic problems, for which the procedures require not only multiplications and reciprocals, but also additions.

### Multiplications and reciprocals

Let us examine problems 1 and 4, which deal with the cost of the digging of a trench. The original language is a mixture of Akkadian, in syllabic notation, and Sumerian logograms. Here, the choice of the language by the ancient scribe seems to have a special meaning in relation to the problem of interest for the present discussion. Thus, in my edition below, translations of Akkadian words are italicized and translations of Sumerian words are not. The sequences between brackets are reconstructions of damaged portions of the text. The numbers on the left correspond to the line numbers in Neugebauer and Sachs' edition.

**YBC 4663 #1**

1. A trench. 5 ninda the length, 1 1/2 ninda <the width>, 1/2 ninda its depth, a volume of 10 <gin> the work assignment, 6 še [of silver the wages].
2. The base, the volume, the (number) of workers and the silver how much? You, in your procedure,
3. the length and the width multiply together, 7.30 *will be given to you*.
4. 7:30 to its depth *raise*, 45 *will be given to you*.
5. The reciprocal of the work assignment loosen, 6 *will be given to you*. *To* 45 *raise*, 4:30 *will be given to you*.
6. 4:30 *to the wages raise*, 9 *will be given to you. Such is the procedure*.

The statement (line 1) gives the dimensions of a trench, the volume of earth that must be extracted every day by a worker ("assigned volume") and the daily wage. It is required to calculate the area of the base of the trench, its volume, and the total wages to be paid to the workers (implicitly: to do the job in one day). The procedure, opened by the formula "You, in your procedure", prescribes a succession of operations: multiply the length by the width to get the base (line 3), then multiply the result by the height to get the volume (line 4), then divide the volume of the trench by the volume that must be extracted by each worker every day to get the number of workers (line 5), and finally, multiply the number of workers by the daily wage to get the total cost (line 6).The final result is 9, and the procedure is closed by the formula "Such is the procedure".

What is of interest for the present discussion is the way in which numbers and quantities are handled. The data of the problem are expressed in specific quantities (measurements of length, volume and weight) in the statement. However, in the procedure itself, only abstract numbers are involved. The relationship between the measurement values given in the

[18] Tablet YBC 4663 was published in [8], pp. 69-75,pl 7 and 32.

statement and the abstract numbers used in the procedure is exactly that established in the metrological tables as shown by Table 5 below[19]:

| | | | |
|---|---|---|---|
| length | 5 *ninda* | → 5 | (table L) |
| width | 1 1/2 *ninda* | → 1:30 | (table L) |
| depth | 1/2 *ninda* | → 6 | (table Lh) |
| wage | 6 *še* | → 2 | (table W) |
| assigned volume | 10 *gin* | → 10 | (table S) |

**Table 5: Data instatement of YBC 4663 #1 and conversions in SPVN with metrological tables**

Then the procedure is a sequence of multiplications with these abstract numbers or their reciprocal (Table 6).

| | | |
|---|---|---|
| length × width | 5 × 1:30 = 7:30 | line 3 |
| surface × depth | 7:30 × 6 = 45 | line 4 |
| volume/ assigned volume | 45 / 10 = 45 × 6 = 4:30 | line 5 |
| number of workers × daily wage | 4:30 × 2 = 9 | line 6 |

**Table 6: Operations in procedure of YBC 4663 #1**

If we observe carefully how the measurement values and the abstract numbers are distributed in the problem, we detect a structure similar to that described in §5 for school exercises. The statement, noted only with Sumerian logograms, provides the data as measurement values. The procedure describes operations with mainly Akkadian terms; these operations are multiplications and reciprocals acting on abstract numbers. The shift from measurement values (in the statement) to abstract numbers (in the procedure) is mirrored by a shift from Sumerian to Akkadian. It seems that the bilingual characteristic to scholarly texts is used here for didactic purposes. Another feature seems to reflect the didactic nature of the problem: Table 5 shows that the first step of the resolution, which is the transformation of measurement values into abstract numbers, involves almost all the metrological tables (only the table for capacity is missing).

Problems 2 to 6 derive from problem 1 by circular permutations of the data. These problems have the same structure as problem 1, but they include an additional step: after the procedure is finished, the abstract number sought is transformed into a measurement value, such that the problem is completely solved. For example, consider problem 4.

### YBC 4663 #4

20. 9 *gin* the (total expenses in) silver for a trench. 5 *ninda* the length, 1 1/2 *ninda* the width. 10 (*gin*) the assigned volume. 6 *še* (silver) the wage (per worker).
21. Its depth how much? In your procedure:
22. The length and the width cross. 7:30 *will be given to you. The reciprocal of the assigned volume detach*,
23. *to* 7:30 *raise*. 45 *will be given to you*. 45 *to the wage raise*.
24. 1.30 *will be given to you. The* reciprocal *of* 1:30 *detach*. 40 *will be given to you*.
25. 40 *to* 9, the (total expense in) silver *raise*. 6, the depth, *will be given to you*. 1/2 *ninda* is its depth.

The structure of the problem and the procedure are similar to that of problem 1. I just summarize the data and operations in Table 7 and Table 8, and will focus on the last additional step.

[19] In what follows, the metrological table for lengths is referred to as "table L", for heights (and other vertical dimensions) as "table Lh", for weights as "table W", for surfaces as "table S", and for capacities as "table C". These tables are available online here ([13], section 9)

| | | | |
|---|---|---|---|
| silver | 9 *gin* | → 9 | (table W) |
| length | 5 *ninda* | → 5 | (table L) |
| width | 1 1/2 *ninda* | → 1:30 | (table L) |
| assigned volume | 10 *gin* | → 10 | (table S) |
| wage | 6 *še* | → 2 | (table W) |

**Table 7: Data in statement of YBC 4663 #4 and conversions in SPVN with metrological tables**

| | | |
|---|---|---|
| length × width | 5 × 1:30 = 7:30 | ligne 22 |
| base/ assigned volume | 7:30 / 10 = 7:30 × 6 = 45 | lignes 22-23 |
| (base / assigned volume) × wage | 45 × 2 = 1:30 | ligne 23-24 |
| silver / [(base / assigned volume) × wage] | 9 / 1:30 = 9 × 40 = 6 | ligne 25 |

**Table 8: Operation in procedure of YBC 4663 #4**

At the end of the procedure, the solution found is an abstract number (6). There, and only there, is the question of the order of magnitude raised. The transformation of the abstract number 6 into a measurement of depth requires the reading of the metrological table for vertical dimensions (Lh). At the same time, a rough evaluation of the order of magnitude of the sought depth is necessary to locate the right portion of the metrological table to be used. The metrological table Lh (see [13], section 9) gives several measurement values associated to the number 6. These possible measurement values are 3 *šu-si* (ca. 5 cm), ½ *ninda* (ca. 3 m), 40 *ninda* (ca. 240 m) and 1 *danna* (ca. 10,5 km). The right choice is clearly ½ *ninda*.

Finally, we see that the process is even closer to that described in §4 than that I suggested just above about problem 1. The spatial disposition which makes the process apparent in the school tablet UM 29-15-192 is replaced here, in YBC 4663, by linguistic markers: the statement is noted in Sumerian, the procedure, noted in Akkadian, is open and closed with a conventional formula, and the answer is noted in Sumerian. These linguistic markers highlight the structure of the process: First the measurement values provided in the statement are transformed into abstract numbers through the direct reading (from left to right) of the metrological tables, second the procedure enumerates the multiplications and reciprocals to be performed with abstract numbers, and finally the result, an abstract number, is transformed into a measurement value through the reverse reading (from right to left) of the metrological tables, with a rough mental evaluation of the orders of magnitude.

## The intervention of addition

The final two problems seem to be similar to the previous one, but a difference must be noted: one piece of information provided by the statement is the sum of the length and width. This detail is of a great consequence: we are dealing with quadratic problems instead of linear problems as in previous examples. The operations in SPVN include not only multiplications and reciprocals, but also additions and subtractions. How are these operations performed with SPVN? We focus here on Problem 7.

## YBC 4663 #7

1. 9 *gin* is the (total expenses in) silver for a trench.
2. The length and the width I added, it is 6:30. ½ *ninda* [its depth].
3. 10 *gin* the assigned volume, 6 *še* (silver) the wage. The length and the width how much?
4. You, for knowing it. The reciprocal of the wage *detach.*
5. *To* 9 *gin*, the silver, *raise.* 4:30 *will be given to you*.
6. 4.30 *to* the assigned volume *raise*. 45 *will be given to you.*
7. The reciprocal of its depth detach. *To* 45 *raise.* 7:30 *will be given to you.*
8. ½ of the length *and* the sag *which I added break.* 3:15 *will be given to you.*
9. 3:15 cross itself. 10:33:45 *will be given to you.*
10. 7:30 *from* 10:33:45 *tear out.*
11. 3:3:45 *will be given to you. Its square root take.*
12. 1:45 *will be given to you. To the one append, from the other cut off.*
13. The length and the width *will be given to you*. 5 (*ninda*) the length, 1 ½ *ninda* the width.

As before, the operations contained in the procedure apply only to abstract numbers, and thus the first step is the transformation of the meteorological data of the statement into SPVN by use of metrological tables (Table 9).

| | | | |
|---|---|---|---|
| silver | 9 *gin* | → 9 | (table W) |
| length + width | | → 6:30 | (provided directly as SPVN) |
| depth | 1/2 *ninda* | → 6 | (table Lh) |
| assigned volume | 10 *gin* | → 10 | (table S) |
| wage | 6 *še* | → 2 | (table W) |

**Table 9: Data in statement of YBC 4663 #7 and conversions in SPVN with metrological tables**

The procedure starts with the calculation of the base, by a method identical to that of the previous problem. In line 7, the base (length × width) is obtained, it is 7:30. The problem is therefore reduced to finding two numbers, corresponding to the length and the width, knowing their sum and product. The data are then the following.

length + width → 6:30
length × width → 7:30

To perform the rest of the procedure, in particular the subtraction in line 10, the number 7:30, which represents the surface, must be placed relative to the number 6:30, which represents the sum of the lengths and the width. Table 10 shows some of the possible relative positions of 7:30 and 6:30; these numerical configurations are designated A, B and C. For each configuration, the place of the units (U in Table 10) must be determined, in order to display the product of two numbers relative to these numbers. This is a new situation in comparison with the previous problems, which involve only multiplications.

| | | **A** | | | **B** | | **C** | | | |
|---|---|---|---|---|---|---|---|---|---|---|
| | | | | **U** | **U** | | **U** | | | |
| 1 *ninda* | → | | 1 | | 1 | | | 1 | | |
| length | → | | 5 | | 5 | | | 5 | | |
| width | → | | 1 | 30 | 1 | 30 | | 1 | 30 | |
| surface | → | 7 | 30 | | 7 | 30 | | | 7 | 30 |
| length + width | → | | 6 | 30 | 6 | 30 | | 6 | 30 | |

**Table 10: some of the possible configurations of numbers in YBC 4663 #7**

Each of these configurations is determined by the choice of the relative positions of the number1, which corresponds to the measurement unit of length 1 *ninda*, and the unit in the number (U), in which I designate in the following as "computing unit" to distinguish it from the measurement units.

The positions are displayed here in columns that evoke a calculation device. This presentation is of course significant to me. It is likely that operations were performed using a physical instrument, and the representation above may give a good idea of how the problems dealing with the relative positions of numbers appeared to the eyes of ancient scribes. These different configurations of the numbers on the "abacus" reflect different possible geometric situations. But nothing in the statement of the problem permits us to favor one position over another.

Let us examine the procedure described in the cuneiform text, and focus on the operations. If we apply, step by step, the calculations prescribed by the procedure in the configuration A (i.e., in 6:30, 30 is placed in the "computing unit" column ), the calculation works well and we obtain the expected value 5, which corresponds to the length 5 *ninda*, and the value 1:30, which corresponds to the width 1 1/2 *ninda*. If we apply, step by step, the calculations prescribed by the procedure in the configurations B (i e., in 6:30, 6 is placed in the "computing unit" column), the calculation works well again, and we also obtain the expected result. It is the same for configuration C, and any configuration (see Table 11).

| | | A | | | B | | | C | | | | |
|---|---|---|---|---|---|---|---|---|---|---|---|---|
| **Lines** | **Operations** | | | **U** | **U** | | | **U** | | | | |
| 7. | Half ( ?) | | 6 | 30 | 6 | 30 | | | 6 | 30 | | |
| | = | | 3 | 15 | 3 | 15 | | | 3 | 15 | | |
| 8. | Product | | 3 | 15 | 3 | 15 | | | 3 | 15 | | |
| | | | 3 | 15 | 3 | 15 | | | 3 | 15 | | |
| | = | 10 | 33 | 45 | 10 | 33 | 45 | | | 10 | 33 | 45 |
| 9. et 10. | Subtraction | 10 | 33 | 45 | 10 | 33 | 45 | | | 10 | 33 | 45 |
| | | **7** | **30** | | **7** | **30** | | | | **7** | **30** | |
| | = | 3 | 3 | 45 | 3 | 3 | 45 | | | 3 | 3 | 45 |
| 10. et 11. | Square root | 3 | 3 | 45 | 3 | 3 | 45 | | | 3 | 3 | 45 |
| | = | | 1 | 45 | 1 | 45 | | | 1 | 45 | | |
| 11. et 12. | Addition | | 3 | 15 | 3 | 15 | | | 3 | 15 | | |
| | | | 1 | 45 | 1 | 45 | | | 1 | 45 | | |
| | = | | 5 | | 5 | | | | 5 | | | |
| | Subtraction | | 3 | 15 | 3 | 15 | | | 3 | 15 | | |
| | | | 1 | 45 | 1 | 45 | | | 1 | 45 | | |
| | = | | 1 | 30 | 1 | 30 | | | 1 | 30 | | |

**Table 11: some of the possible configurations of the numbers in the procedure in YBC 4663 #7**

Some conclusions can be drawn from this example.
1) There is no requirement to choose that the measurement unit 1 *ninda* is represented on the abacus by the computing unit (1 placed in column U). We can decouple the measurement unit *ninda* and the computing unit.
2) Execution of the calculation requires only the initial determination of the positions on the “abacus” of some minimum data, for example the relative positions of the computing unit and the number 1 representing 1 *ninda*.
3) The problem of the order of magnitude only appears at the end of the calculation, when the result (here, 5) should be interpreted as a measure of length (here, 5 *ninda* - see line 12).

This example shows that the performance of the calculation does not require the tracking of orders of magnitude of the numbers placed on the abacus, but only their relative positions. At the same time, these numbers represent quantities. Numbers and operations have a geometric meaning, and these meanings can guide the reasoning that governs the calculation. Each step of the calculation can thus be represented by operations on squares and rectangles (see Figure 7).[20]

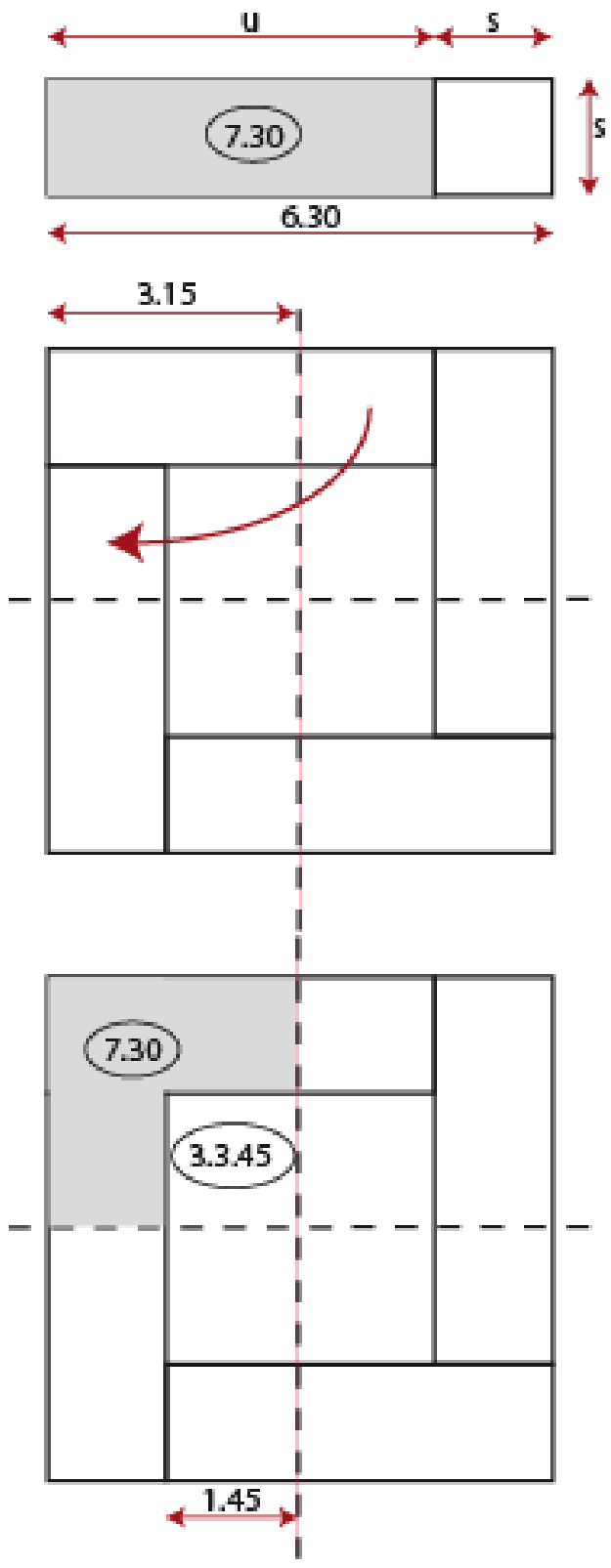

**Figure 7: geometrical representation of the procedure in YBC 4663 #7**

## 7- Orders of magnitude and place of the unit in the number

The examples described in this brief presentation show that the floating positional sexagesimal numbers used in Old Babylonian period were basically devices for calculating, not for representing quantities. In the ideal world of school mathematics, only multiplications and reciprocal are performed with these floating numbers. Moreover, numbers in SPVN are usually regular. The power and efficiency of the fundamental factorization algorithms, which allow the calculation of reciprocals, square roots, and cube roots, rely on the floating nature of the numbers on which they are implemented. Many of the problems encountered in cuneiform mathematics, including in scholarly contexts, are linear and are solved by a series of multiplications and reciprocals. Some mathematical methods used in the Old Babylonian period are particularly well suited to this multiplicative universe: false position, enlargement or reduction of figures, and the use of coefficients.

However, the multiplicative universe is too limited to solve some problems, including quadratic problems where additions and subtractions intervene. For the needs of certain calculations, the position of the computing units must be fixed at specific times of the calculation. Fixing the position of computing units is equivalent to setting rules for the use of the abacus. These rules, which reflect geometric situations, leave much more room for flexibility in calculations than modern notation with marks such as 0, semi colon (;), colon (,)

[20] The fact that the quadratic procedures are based on operations on squares and rectangles was shown by Jens Høyrup – see for example [6].

or the extended degree-second-minute system (``,`, °,′, ″ ) currently used in modern publications.

In our decimal number system, and its associated decimal metric system, the orders of magnitude of a quantity are fixed by defining the position of the unit in the number which measures the quantity. But this is not the case in the Mesopotamian scholarly system that we have just described. Determining the order of magnitude is an aspect of metrology. Placing the numbers relative to each other is an aspect of the calculation itself, which was probably carried out on an abacus

# Appendix: numbers and measurement units

## 1- Numbers

### a- System S (sexagesimal and additive)

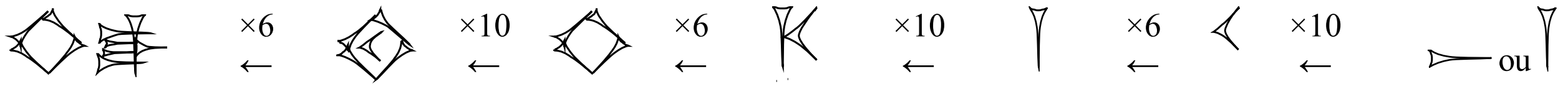

| *šargal* | *šaru* | *šar* | *gešu* | *geš* | *u* | *Aš* or *diš* |
|---|---|---|---|---|---|---|
| 216000 | 36000 | 3600 | 600 | 60 | 10 | 1 |

### b- Common fractions

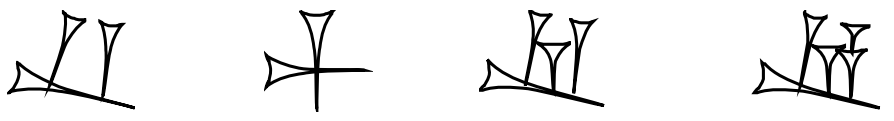

1/3 1/2 2/3 5/6

### c- Abstract numbers (sexagesimal place value notation = SPVN)

Signs

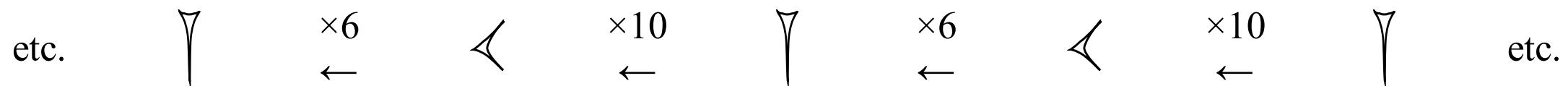

Sexagesimal place value notation in cuneiform texts uses two signs: the wedge (1) and the oblique wedge (10) - see above. The sexagesimal digits from 1 to 59 are written by juxtaposition of 1 and 10 as many times as is necessary (additive decimal notation). Example:

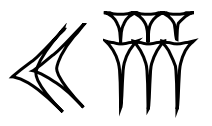 = 26

In a number with several sexagesimal positions, a sign written in a position is sixty times more than the same sign written in the previous position (i.e. located at its right).

There is no notation in cuneiform texts equivalent to our "zero" or the mark we use to separate the integer part from the fractional part in a number (equivalent to our decimal point, as in "3.14"). In other words, cuneiform notation does not indicate the position of the units in numbers. For example, the same sign 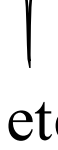 represents the numbers we write 1, 60, 3600, 1/60, etc. in base ten, or 1, 1:0, 1:0:0 or 0;1 etc. in base sixty.

For example, a relationship between floating numbers such as 2×30 = 1expresses the fact that30is the reciprocal of 2 (their product is "1"). Similarly 44:26:40×1:21 = 1 means that44:26:40is the reciprocal of1:21. The relationship 30×30 = 15 expresses the fact that the (floating) square of 30 is 15, or the (floating) square root of 15 is 30. Strictly speaking, we should replace the sign '=' by a sign for modular equality '≡', it being understood that in the equivalence relation in question here, two numbers are equivalent if their quotient is a power of 60, positive or negative integer exponent.

## 2- Measurement units

The following diagrams are a synthetic representation (in part following Jöran Friberg) of information provided by the metrological tables. They represent the factors between measurement units, the graphemes representing these measurement units, and the corresponding abstract numbers.

**Capacity** (1 *sila*≈1 litre)

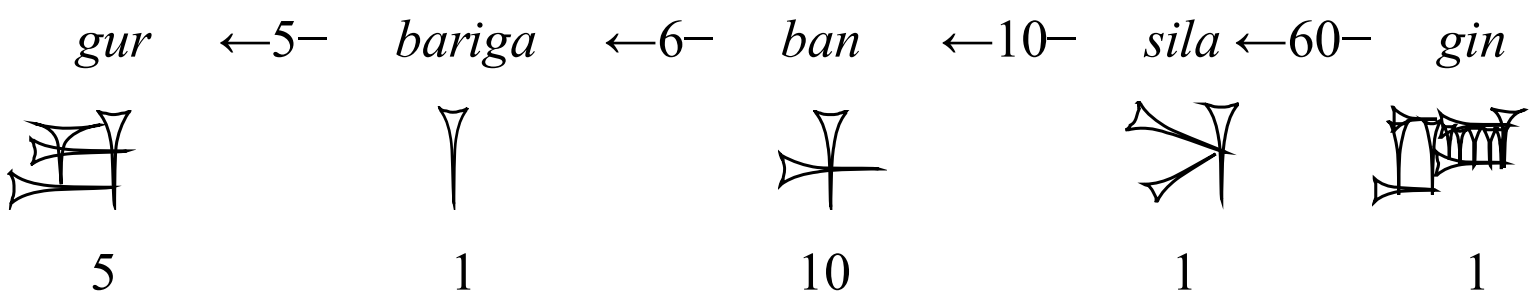

**Weight** (1 *gu*≈30 kg)

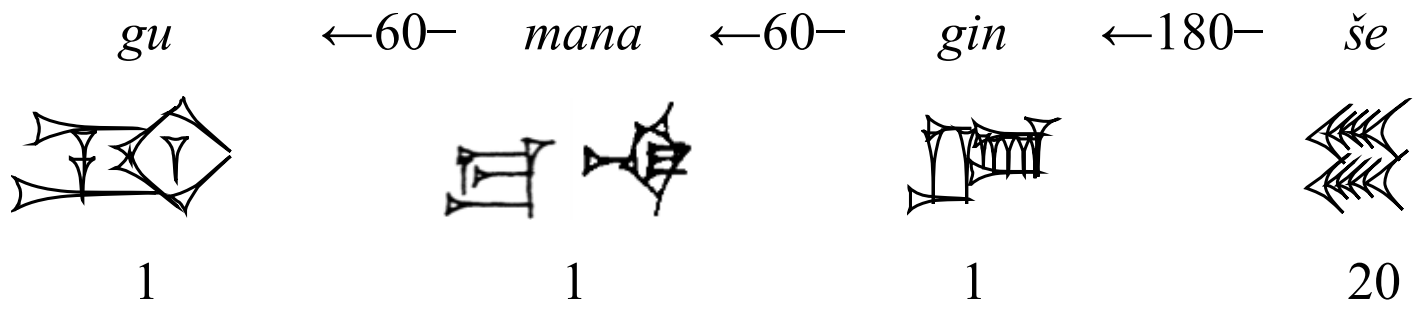

**Surfaces** (1 *sar*≈36 m ) **and volumes** (1 *sar*-volume = 1 *sar* with thickness1 *kuš* ≈ 18 $m^3$)

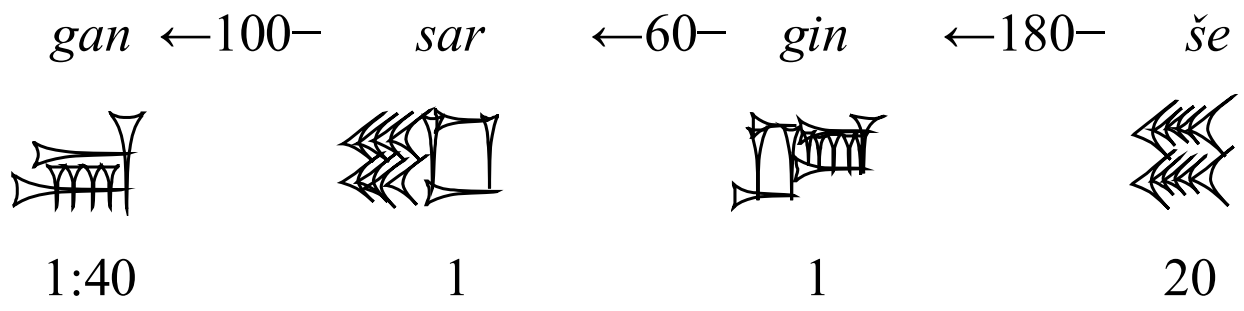

**Lengths** (1 *ninda*≈6 m)

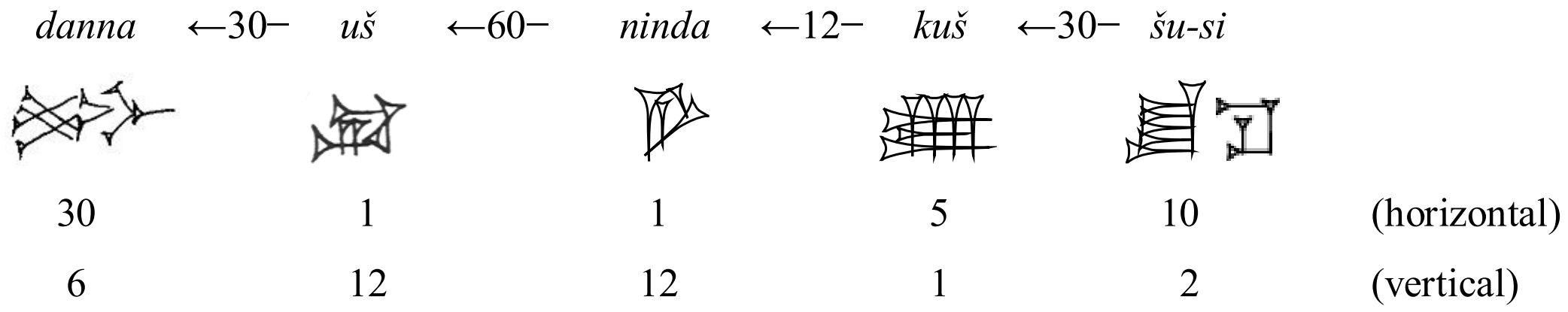

## 3- *MesoCalc*, a Mesopotamian calculator

The representations of numbers and quantities reflect how they are used in calculations. In cuneiform documents, the notations have been shaped by ancient scribes for specific uses; this article is an attempt to describe these uses. The conventions for the representation of the cuneiform signs adopted by modern historians reflect how they understand the ancient calculations. The elaboration of these conventions is therefore of great importance, and is often underestimated in historical studies.

Working with the cuneiform mathematical texts requires great computing power that modern tools allow us to model with computer programs. For the reasons stated above, the choice of notations and calculation principles adopted with computer programs rely on historical considerations. MesoCalc, the tool developed by Baptiste Mélès (http://baptiste.meles.free.fr/site/mesocalc.html – see Appendix 3) is an attempt to implement and automate calculations as faithfully as possible to ancient practices. In particular, this tool distinguishes metrological calculations and calculations in sexagesimal place value notation, and implements the floating calculation for multiplication and reciprocals.

## Sources

| Museum number | Publication | CDLI number and link to CDLI |
|---|---|---|
| CBS 1215 | [18], 238-240 | P254479 |
| HS 241 | [5], n°42 ; [12]; n°29 | P388160 |
| Ist Ni 10241 | [10] | P368962 |
| MS 3874 | [4] | P252955 |
| UM 29-15-192 | [9], 248, 251 | P254900 |
| YBC 4663 | [8], 69 (text H) | P254984 |
| YBC 7302 | [8], 44 | P255051 |